\documentclass[11pt]{amsart}
\usepackage{color}
\usepackage{epsfig}
\usepackage{psfrag, graphics}
\usepackage{amsthm}
\usepackage{amssymb}
\usepackage{amsmath}
\usepackage{amscd}
\usepackage{graphicx}
\usepackage{epstopdf}
\newcommand{\eqdef}{\;{:=}\;}

\newtheorem {Theorem}   {Theorem}

\numberwithin{Theorem}{section}

\newtheorem {Lemma}[Theorem]    {Lemma}
\newtheorem {Proposition}[Theorem]{Proposition}
\newtheorem {Corollary}[Theorem]{Corollary}
\newtheorem {Question}[Theorem]{Question}
\theoremstyle{definition}

\theoremstyle{remark}
\newtheorem{Remark}[Theorem]{Remark}
\newtheorem{Example}[Theorem]{Example}

\expandafter\chardef\csname pre amssym.def
at\endcsname=\the\catcode`\@ \catcode`\@=11
\def\undefine#1{\let#1\undefined}
\def\newsymbol#1#2#3#4#5{\let\next@\relax
 \ifnum#2=\@ne\let\next@\msafam@\else
 \ifnum#2=\tw@\let\next@\msbfam@\fi\fi
 \mathchardef#1="#3\next@#4#5}
\def\mathhexbox@#1#2#3{\relax
 \ifmmode\mathpalette{}{\m@th\mathchar"#1#2#3}%
 \else\leavevmode\hbox{$\m@th\mathchar"#1#2#3$}\fi}
\def\hexnumber@#1{\ifcase#1 0\or 1\or 2\or 3\or 4\or 5\or 6\or 7\or 8\or
 9\or A\or B\or C\or D\or E\or F\fi}

\font\teneufm=eufm10 \font\seveneufm=eufm7 \font\fiveeufm=eufm5
\newfam\eufmfam
\textfont\eufmfam=\teneufm \scriptfont\eufmfam=\seveneufm
\scriptscriptfont\eufmfam=\fiveeufm

\catcode`\@=\csname pre amssym.def at\endcsname

\newcommand{\AC}{{\mathcal A}}

\newcommand{\FF}{{\mathcal F}}
\newcommand{\GG}{{\mathcal G}}
\newcommand{\HH}{{\mathcal H}}

\newcommand{\JJ}{{\mathcal J}}

\newcommand{\LL}{{\mathcal L}}
\newcommand{\MM}{{\mathcal M}}

\newcommand{\PP}{{\mathcal P}}

\newcommand{\RR}{{\mathcal R}}

\def    \C      {{\mathbb C}}
\def    \R      {{\mathbb R}}
\def    \Z      {{\mathbb Z}}

\def    \Q      {{\mathbb Q}}

\newcommand{\id}{{\operatorname{id}}}

\def    \om       {\omega}
\def    \eps      {\epsilon}
\def    \ssminus  {\smallsetminus}
\def    \p        {\partial}

\def    \rk       {\mathit{rk}\,}

\def    \SB      {\operatorname{SB}}

\def    \Hh     {\mathbf{H}}

\def    \12      {\frac{1}{2}}

\def    \sl        { \Phi_{L}} 
\def    \sr         { \Phi_{R}}
\def    \sll         { \Phi_{L,k}} 
\def    \srl        { \Phi_{R,k}} 
\def    \ind      {\operatorname{ind}}
\def    \PHG         { \Psi_{\scriptscriptstyle{HG}}} 
\def  \phg              {\varrho_t^{\scriptscriptstyle{HG}}}
\def  \phg              {\varrho_t^{\scriptscriptstyle{HG}}}

\def  \pfg              {\varrho_t^{\scriptscriptstyle{FG}}}

\def    \H   {\operatorname{H}}  
\def    \CF   {\operatorname{CF}} 
\def    \CM   {\operatorname{CM}} 
\def    \CZ   {\operatorname{\mu_{\scriptscriptstyle{CZ}}}} 

\begin{document}







\title[Hofer's geometry and Floer theory under the quantum limit]{Hofer's geometry and Floer theory under the quantum limit}

\author[Ely Kerman]{Ely Kerman}
\address{Department of Mathematics, University of Illinois at Urbana-Champaign, Urbana, IL 61801, USA }
\email{ekerman@math.uiuc.edu}

\date{\today}

\thanks{This research was partially supported by NSF Grant DMS-0405994 and a grant from the Campus Research
Board of the University of Illinois at Urbana-Champaign.}


\bigskip

\begin{abstract}
In this paper, we use Floer theory to study the Hofer length functional 
for paths of Hamiltonian diffeomorphisms which are sufficiently short.
In particular, the length minimizing 
properties of a short Hamiltonian path are related to the properties and number of its 
periodic orbits.    
\end{abstract}

\maketitle

\section{Introduction}

Throughout this work, $(M,\om)$ will be a closed symplectic manifold of dimension $2n$.
The space of smooth $\om$-compatible almost complex structures on $M$ will be denoted by $\JJ(M,\om)$. 
For $J$ in $\JJ(M,\om)$, let $\hbar(J)$ be the infimum\footnote{We use the convention that the infimum over the empty set is equal to infinity.} 
over the symplectic areas of nonconstant $J$-holomorphic 
spheres in $M$, and set
$$
\hbar = \sup_{J \in \JJ(M,\om)} \hbar(J).
$$
This strictly positive quantity is the \emph{quantum limit} referred to in the title. In this paper, we consider  
paths of Hamiltonian diffeomorphisms whose Hofer length is less than $\hbar$, and prove several results 
which relate the length minimizing properties of these paths to the periodic orbits of the
corresponding Hamiltonian flows. To establish these results we develop some new Floer theoretic tools. 
These are motivated, in part, by similar constructions in Lagrangian Floer theory due to Chekanov, \cite{ch}, 
as well as recent work by Albers in \cite{al}.
These tools allow one to detect length minimizing  
Hamiltonian paths using homological algebra, \cite{kl,oh2,oh3,schw}.

\subsection{Definitions}

Before stating the main results, we first recall some basic definitions and fix some conventions.
For more detailed accounts of this material the reader is referred to the books 
\cite{hz:book,mcsa:book,po:book}.  

Let $S^1 = \R /\Z$  be the circle parameterized by $t \in [0,1]$. A
{\bf Hamiltonian} on $M$ is a smooth function $H$ in $C^{\infty}(S^1 \times M)$, which we
also view as a loop of smooth functions $H_t(\cdot)=H(t,\cdot)$ on $M$. We say that $H$ is 
normalized if $\int_M H_t \, \om^n=0$ for each $t \in [0,1].$ The space of 
normalized Hamiltonians is denoted by $C^{\infty}_0(S^1 \times M)$.

Each Hamiltonian $H$ determines a time-dependent vector
field $X_H$ on $M$ via the equation
$$
\om(X_H, \cdot) = -dH_t(\cdot).
$$
The flow of this vector field (also called the Hamiltonian flow of $H$) is denoted by
$\phi^t_H$ and is defined for all $t \in [0,1]$. The group of Hamiltonian diffeomorphisms
of $(M,\om)$ consists of all the time one maps of such flows.

For every path of Hamiltonian diffeomorphisms, $\psi_t$, there is a unique normalized Hamiltonian 
$H \in C^{\infty}_0(S^1 \times M)$ such that $\phi^t_H \circ \psi_0 = \psi_t.$ Following Hofer \cite{ho1}, this generating Hamiltonian
is used to define the Hofer length of the 
path $\psi_t$ as follows
\begin{eqnarray*}
\operatorname{length} (\psi_t)&=& \|H\|\\
{} &=& \int_0^1 \max_M H_t \,\,dt-  \int_0^1 \min_M H_t \,\,dt \\
{} &=& \|H\|^+ + \|H\|^-.
\end{eqnarray*}
Both $\|H\|^+$ and  $\|H\|^-$ provide different measures of the 
length of $\psi_t$, called the positive and negative Hofer lengths, respectively.

Let $[\psi_t]$ be the class of Hamiltonian paths which are homotopic to $\psi_t$ relative to its endpoints. Denote the set of normalized Hamiltonians 
which generate the paths in $[\psi_t]$ by 
$$
C^{\infty}_0([\psi_t]) =\{ H \in C^{\infty}_0(S^1 \times M) \mid [\phi^t_H \circ \psi_0] = [\psi_t]\}.
$$  
The Hofer seminorm of $[\psi_t]$ is then defined by 
$$
\rho_{\scriptscriptstyle{\H}}([\psi_t]) = \inf \{ \|H\| \mid H \in C^{\infty}_0([\psi_t]) \}.
$$
The positive and negative Hofer seminorms of $[\psi_t]$ are defined similarly as
$$
\rho^{\pm}([\psi_t]) = \inf\{ \|H\|^{\pm}\mid H \in C^{\infty}_0([\psi_t]) \},
$$
and we call
$$
\bar{\rho}([\psi_t]) = \rho^+([\psi_t]) + \rho^-([\psi_t])
$$
the two-sided Hofer seminorm of $[\psi_t]$.

Since these seminorms are bi-invariant, we need only consider paths which start at the identity
and hence have the form $\phi^t_H$ for some $H$ in $C^{\infty}_0(S^1 \times M)$. To avoid 
pathologies, we will also assume that our paths are {\bf regular}.  That is, we will only consider 
Hamiltonian paths $\phi^t_H$ for which the functions $H_t$ do not vanish identically
for $t \in (0,1)$.

We say that $\phi^t_H$ minimizes the Hofer length in its homotopy class
if there is no path in $[\phi^t_H]$ with a smaller  Hofer length, i.e., $\|H\|= \rho_{\scriptscriptstyle{\H}}([\phi^t_H])$. 
The notion of a path which minimizes the positive, negative or two-sided Hofer length is defined
analogously. Note that $\rho_{\scriptscriptstyle{\H}}([\psi_t]) \geq \bar{\rho}([\psi_t])$ and so if  $\phi^t_H$ minimizes $\bar{\rho}$ then it also minimizes $\rho_{\scriptscriptstyle{\H}}$. 
It should also be mentioned that length 
minimizing paths for the Hofer seminorm need not exist in a given homotopy class.
Explicit examples of such classes are constructed by Lalonde and McDuff in \cite{lm}.

We will use Floer theory to relate the length minimizing properties of a Hamiltonian path to the number and 
properties of its periodic orbits. Hence, we will focus our attention on $\PP(H)$, the set of contractible periodic orbits of $X_H$ which are $1$-periodic (have period equal to one). We say that $x$ in $\PP(H)$
is nondegenerate if the map $d\phi^1_H \colon T_{x(0)}M \to T_{x(0)}M$ does not have 
one as an eigenvalue. A Hamiltonian $H$ will be called a 
{\bf Floer Hamiltonian} if each $x$ in $\PP(H)$ is nondegenerate.

A spanning disc for a $1$-periodic orbit $x$ is a smooth map $u$ 
from the unit disc $D^2 \subset \C$ to $M$, such that 
$u(e^{2\pi i t}) = x(t)$. It can be used to define  
two important quantities for $x$ which can be studied using Floer theory. The first  quantity  is $\CZ(x,u)$, 
the Conley-Zehnder index of $x$ with respect to $u$. This index is normalized here so that if $p$ is a nondegenerate critical point 
of a $C^2$-small time-independent Hamiltonian, and $u$ is the constant spanning disc, then $$\CZ(p,u)=\ind(p)-n,$$ where 
$\ind(p)$ is the Morse index of $p$.
The second important quantity is the action of $x$ with respect to $u$, which
is defined by
$$
\AC_H(x,u) = \int_0^1 H(t,x(t))\,dt - \int_{D^2} u^* \om.
$$
Both the Conley-Zehnder index and the action of $x$ with respect to $u$ depend only on the homotopy class of $u$.


\subsection{Results} 

If $\phi^t_H$ does not minimize $\rho_{\scriptscriptstyle{H}}$ in its homotopy class,
then it also fails to minimize  $\rho^{+}$ or $\rho^-$. For this reason, we
formulate and prove our results for the one-sided seminorms. In particular, we consider 
the positive Hofer seminorm. The statements and proofs of the corresponding results for $\rho^-$ are entirely similar
(see Corollary \ref{flip}).

\begin{Theorem}\label{general}
Let $H$ be a Floer Hamiltonian such that $\|H\| < \hbar$. 
If $\phi^t_H$ does not minimize the positive Hofer seminorm in its 
homotopy class, then there are at least $\rk(\H(M;\Z))$ contractible $1$-periodic orbits $x_j$ of $H$ 
which  admit spanning disks $u_j$ such that 
$$ -n \leq\CZ(x_j,u_j) \leq n$$
and  
$$-\|H\|^- \leq \AC_H(x_j,u_j) < \|H\|^+.$$ 
\end{Theorem}
\noindent Here, $\rk(\H(M;\Z))$ is the rank of the abelian group $\H(M;\Z)$, and the rank of the torsion
subgroup is the minimal number of elements needed to generate it. 

For a general, possibly degenerate, Hamiltonian we prove: 
\begin{Theorem}\label{degenerate}
Let $H$ be a Hamiltonian such that $\|H\| < \hbar$. If $\phi^t_H$ does not minimize
the positive Hofer length in its homotopy class, then there is a
 contractible $1$-periodic orbits $y$ of $H$ 
which  admits  a spanning disk $w$ such that 
$$-\|H\|^- \leq \AC_H(y,w) < \|H\|^+.$$
\end{Theorem}

In the corresponding results for the negative Hofer seminorm, the action values
are confined  instead to the interval $(-\|H\|^-, \|H\|^+]$, (see, again, Corollary \ref{flip}).

To refine Theorems \ref{general} and \ref{degenerate} it is useful to restrict ones attention to Hamiltonians which have properties
known to be necessary to generate a length minimizing path.
Following \cite{bp}, a Hamiltonian $H$ is called {\bf quasi-autonomous} if it has at least one fixed 
global maximum $P \in M$ and one fixed global minimum $Q \in M$. That is, 
$$
H(t,P) \geq H(t,p) \geq H(t,Q)
$$
for all  $p \in M$ and $t \in [0,1]$.


\begin{Theorem}\cite{bp,lm}
If the  Hamiltonian path generated by $H$ minimizes the Hofer norm in its homotopy class, then $H$ must be 
quasi-autonomous.  
\end{Theorem}

A symplectic manifold $(M,\om)$ is said to be {\bf spherically rational}, if the quantity
$$
r(M,\om) = \inf_{A \in \pi_2(M)} \left\{ |\om(A)| \mid |\om(A)|>0\right\}.
$$
is strictly positive. In this case, $r(M,\om) \leq \hbar$ and  we prove:

\begin{Theorem}\label{quasi-autonomous}
Let $H$ be a  Floer Hamiltonian on a spherically rational symplectic manifold $(M,\om)$ such that $H$ is quasi-autonomous and $\|H\| < r(M, \om)$.
If $\phi^t_H$ does not minimize $\rho_{\scriptscriptstyle{\H}}$ in its homotopy class,  then $H$ has at least $\SB(M)+2$ 
contractible $1$-periodic orbits. 
\end{Theorem}

Here, $\SB(M)$ is equal to $\rk (H(M;\Q)$, the sum of Betti numbers of $M$. 
The Arnold conjecture for compact symplectic manifolds, \cite{cz,fl1,fl2,hs,fo,lt}, implies the 
existence of at least $\SB(M)$ contractible $1$-periodic orbits of a Floer Hamiltonian. The two \emph{extra} orbits detected in Theorem \ref{quasi-autonomous}
were also found in the case of symplectically aspherical manifolds  in \cite{kl}. One is lead by these results (as well 
those mentioned below) to 
the following question.
\begin{Question}
If the path generated by a  Floer Hamiltonian does not minimize the (positive/negative) Hofer length (in its homotopy class), must it have at least $\SB(M)+2$ 
contractible $1$-periodic orbits?
\end{Question}
For paths generated by autonomous Hamiltonians there is a well known conjecture in this direction which is motivated by the work of Hofer in \cite{ho2}. It states that if a path generated by an autonomous Hamiltonian is not length minimizing in its homotopy class, then there must be a nonconstant contractible periodic orbit of the Hamiltonian flow with period less than or equal to one, \cite{po:book}.\footnote{It seems reasonable to refine this conjecture to orbits of period equal to one.} This conjecture was essentially settled for all compact symplectic manifolds by McDuff and Slimowitz, in \cite{mcsl}. An extra hypothesis in the main theorem of \cite{mcsl}, concerning the under-twistedness of all critical points, was later shown to be unnecessary by Schlenk in \cite{schl}.   

\subsection{Examples}

Theorems \ref{general}, \ref{degenerate} and \ref{quasi-autonomous} can each be viewed from two perspectives.
In terms of Hofer's geometry, they can be phrased as sufficient conditions for a Hamiltonian path to be  length minimizing for some Hofer seminorm. In terms of dynamics, these results can also be used to detect the existence of special $1$-periodic orbits of a Hamiltonian path which is known not to minimize a Hofer seminorm.  In this section, we describe some elementary examples of Hamiltonian paths on the two sphere which one can show are length minimizing in their homotopy class using the results proved here. We also describe a basic setting wherein 
Theorems \ref{general}, \ref{degenerate} and \ref{quasi-autonomous} can be used to detect periodic orbits for Hamiltonian flows of independent interest.  
 
\subsubsection{Length minimizing Hamiltonian paths on $S^2$}

Consider the unit sphere $S^2$ in $\R^3$. One has coordinates
$$\theta \colon S^2 \ssminus \{p_N,p_S\} \to [0,2\pi]$$  $\text{  and  }$ $$z \colon S^2 \to [-1,1],$$ 
where $p_N$ and $p_S$ are the north and south poles, $\theta(x,y,z)$ is the 
angle of the point $(x,y)$ with respect to the positive $x$-axis, and $z$ is the height.
The area form $\om$ on $S^2$ inherited from the Euclidean volume form on $\R^3$, is a symplectic form 
for which $\hbar= 4\pi$.
 
 The Hamiltonian flow of the function $H_{\upsilon}(\theta,z)=\upsilon z$ is given by 
 $$\phi^t_{H_{\upsilon}}(\theta,z)=(\theta+\upsilon t, z).$$ 
 For $\upsilon < 2\pi$, this flow has no nonconstant periodic 
 orbits of period less than or equal to one.
 It then follows from \cite{mcsl}, that $\phi^t_{H_\upsilon}$
 minimizes the two-sided Hofer length in its homotopy class.  
 This is also implied by Theorem \ref{general}. In particular, the only 
 $1$-periodic orbits of the flow are the constant orbits at the poles, 
 The action of $p_N$ with respect to any spanning disc has the form 
 $\upsilon + 4\pi k$  for some integer $k$. Hence, for $\upsilon < 2 \pi$ the north pole 
 does not admit a spanning disc for which the 
 corresponding action lies in the interval $[-\|H_{\upsilon}\|^-,\|H_{\upsilon}\|^+)=[-\upsilon,+\upsilon).$ Since  $\rk(\H(S^2;\Z))=2$ and only the south pole admits a spanning disk yielding an action value in $[-\upsilon, +\upsilon)$, it follows from Theorem \ref{general} that $\phi^t_{H_\upsilon}$
 minimizes the positive Hofer length in its homotopy class. The analogous result for the negative Hofer seminorm
 implies that $\phi^t_{H_\upsilon}$  also minimizes the negative Hofer length in its homotopy class.

In this setting, one can also construct simple examples of time-dependent 
Hamiltonians which generate length minimizing paths. Consider a Hamiltonian 
$H = H(t,z)$  which does not depend on $\theta$ and has a fixed global maximum (resp. minimum) at 
$z=1$  (resp. $ z=-1$). 
The time-$t$ flow of $H$ is given by $$\phi^t_H(\theta, z)= (\theta + \int_0^t \p_z H(\tau, z) \,d \tau, z).$$ 
If $$0< \int_0^1 \p_z H(t, z) \,d t<2\pi  \text{  for all  } z \in [-1,1],$$ then $\|H\| < 4\pi$ and  the only $1$-periodic orbits of $H$ are again the constant orbits at the poles. Arguing as above, Theorem \ref{general} implies that $\phi^t_H$ minimizes both $\rho^+$ and $\rho^-$, and hence $\bar{\rho}$, in its homotopy class.  Note that these length minimizing paths may have 
many nonconstant periodic orbits with period less than one. They also come in infinite dimensional families (with fixed endpoints).

\begin{Remark}
These examples also illustrate the limitations of the methods developed here.
For example, for $\pi <\upsilon <2\pi$, the path generated by  
$H_{\upsilon} \colon S^2 \to \R$ is not length minimizing among 
all paths (e.g. the flow of $-H_{2\pi - \upsilon}$ has the same time one map). However, there 
is no discernable change in the Floer complexes of the $H_{\upsilon}$ as the parameter 
$\upsilon$ crosses the value $\pi$.
\end{Remark}

\begin{Remark} 
 At least for $S^2$, the problem of characterizing  length minimizing Hamiltonian paths   
above the quantum limit may not be a meaningful.
In particular, it not clear to the author  whether there exists a Hamiltonian $H$ on $(S^2, \om)$ 
such that  $\|H\|>\hbar$ and  $\phi^t_H$ minimizes the Hofer length in its homotopy class. 
\end{Remark}

\subsubsection{Examples of short paths which are not length minimizing}

We now provide some examples of Hamiltonian paths which satisfy the hypotheses 
of the theorems proved in this paper. 
As described below, similar flows appear in several
applications of Hofer's geometry to Hamiltonian dynamics and symplectic topology.

The examples are constructed using 
Sikorav's curve shortening procedure for Hamiltonian paths supported on sets with 
finite displacement energy. 
The displacement energy of a subset $V \subset M$ is defined as
$$
e(V) = \inf_{H \in C^{\infty}_0(S^1 \times M)} \{ \|H\| \mid  \phi^1_H(V) \cap V = \emptyset \}.
$$

\begin{Proposition}\label{shorten}
Let $H$ be an autonomous normalized Hamiltonian which is constant on the complement of an open subset $U$ of $M$ which has finite displacement energy. If $\|H\| > 4 e(U) $, then $$\|H\| > \rho_{\scriptscriptstyle{\H}}([\phi^t_H]).$$ If  $H$ is equal to its minimum value outside of $U$ and $\|H\|^+ > 2 e(U)$, then $$\|H\|^+ \geq \rho^+([\phi^t_H])  + \12 \|H\|^-.$$ 

\end{Proposition}
In the first case, $\phi^t_H$ does not minimize the Hofer length in its homotopy class
and in the second case it does not minimize the positive Hofer length in its homotopy class.
The first part of Proposition \ref{shorten} is proved in \cite{schl} as Proposition 2.1. 
The second part can be obtained by repeating the argument from
\cite{schl} with only minor changes.

Let $N \subset M$ be a closed submanifold of $M$. If $2e(N)<\hbar$, then for sufficiently small
tubular neighborhoods $U$ of $N$ one can easily construct an autonomous normalized  function $H$ which is equal to its minimum value outside of $U$ and satisfies
$\|H\|< \hbar $ and $\|H\|^+>2e(U)$. By Proposition \ref{shorten},  the corresponding path $\phi^t_H$
does not minimize the positive Hofer length. Perturbing this function if necessary, one can 
then apply Theorem \ref{general}, \ref{degenerate}, or
\ref{quasi-autonomous} to obtain information about the periodic orbits of $H$. If $H$ is chosen to be
radially symmetric in the normal directions to $N$,
 then these  periodic orbits are often 
of considerable interest. For example, when $N$ is a Lagrangian submanifold, the flow of $H$ is a reparameterization of a geodesic flow on $N$, 
and when $N$ is a symplectic submanifold, the flow of $H$ can model the motion of a charged particle in a nondegenerate magnetic 
field on $N$. The periodic orbits of such Hamiltonian flows have been studied in several works on Hamiltonian dynamics, Lagrangian embeddings 
and symplectic intersection phenomena, 
(see, for example, \cite{bps,cgk,gi:coisotropic,gu,ke, schl,the,vi:torus}).   
In subsequent work, \cite{ke:coisotropic}, we will consider applications of the techniques developed here to such problems. 

\subsection{Paths with close endpoints}

The methods developed in this work also apply to Hamiltonian paths 
which may be long but whose endpoints are close. The following result
implies Theorem \ref{general} but will be proved separately.

\begin{Theorem}\label{technical}
Let $H$ be a Floer Hamiltonian such that $\bar{\rho}([\phi^t_H]) < \hbar$. 
For every $\epsilon^+>0$ and $\epsilon^->0$, there are at least $\rk(\H(M;\Z))$ 
contractible $1$-periodic orbits $x_j$ of $H$ 
which  admit spanning disks $u_j$ such that 
$$ -n \leq\CZ(x_j,u_j) \leq n$$
and  
$$-\rho^-([\phi^t_H])-\epsilon^- \leq \AC_H(x_j,u_j) \leq \rho^+([\phi^t_H])+ \epsilon^+.$$ 
\end{Theorem}

\begin{Remark}
If $\rho^{+}[\phi^t_H]$  (resp. $\rho^{-}[\phi^t_H]$) is realized by a Hamiltonian path in class $[\phi^t_H]$, then we can set $\epsilon^+=0$ (resp. $\epsilon^- = 0$). One can also set $\epsilon^{\pm}=0$ 
if the symplectic manifold $(M,\om)$ is rational.
\end{Remark}

\subsection{Acknowledgments.}
The author would like to thank Peter Albers for many helpful discussions and for pointing out a gap 
in a previous version of this paper. He would also like to thank
Viktor Ginzburg  and the referee for their insightful comments and suggestions.

\section{Floer theory under the quantum limit}\label{floer}

\subsection{Starting data}

Let $H$ be a Floer Hamiltonian.  Denote the space of smooth $S^1$-families of $\om$-compatible almost complex structures on $M$  by  $\JJ_{S^1}(M,\om)$.
Fixing  a family $J$ in $\JJ_{S^1}(M, \om)$, we refer to the pair $(H,J)$ as our {\bf Hamiltonian data}.

Let $F_s$ be a smooth $\R$-family of functions in $C^{\infty}(S^1 \times M)$ or elements of $\JJ_{S^1}(M,\om)$.
 We call $F_s$ a {\bf compact homotopy} from $F^-$ to $F^+$, if there is a $\tau>0$ such that 
$$
F_s =
\left\{
  \begin{array}{ll}
    F^- ,& \hbox{ for $s \leq  -\tau$ ;} \\
    F^+, & \hbox{ for $s \geq  \tau$ .}
  \end{array}
\right.
$$ 
A {\bf homotopy triple} for our Hamiltonian data $(H,J)$ is a collection of compact homotopies 
$$\HH = (H_s,K_s, J_s).$$ Until we state otherwise,  $H_s$ will be a 
compact homotopy from the zero function to $H$, 
$K_s$ will be  a compact homotopy from the zero function to itself, and $J_s$ will be a compact homotopy  in $\JJ_{S^1}(M,\om)$ from some
$J^-$ to $J$. 

The {\bf curvature} of a homotopy triple $\HH=(H_s,K_s,J_s)$ is the function 
 $\kappa(\HH) \colon \R \times S^1 \times M \to \R$ given by
$$\kappa(\HH) = \p_sH_s - \p_tK_s + \{H_s,K_s\}.$$
We define the positive and negative norms of the curvature by 
$$
|||\kappa(\HH)|||^+ = \int_{\R \times S^1} \max_{p\in
M}\kappa(\HH)\,\,ds \,dt,
$$
and 
$$
|||\kappa(\HH)|||^- = -\int_{\R \times S^1}  \min_{p \in
M}\kappa(\HH) \,\,ds\,dt.
$$
\begin{Example}\label{ex:linear}
Let $\eta \colon \R \to [0,1]$ be a smooth nondecreasing function 
such that $\eta(s) =0 $ for $s \leq -1$ and $\eta(s)=1$ for $s \geq 1$.
A {\bf linear homotopy triple} for $(H,J)$ is a triple of the form 
$$
\overline{\HH}=\left(\eta(s)H, 0, J_s \right).
$$
For a linear homotopy triple we have
$$
\kappa(\overline{\HH})= \dot{\eta}(s)H,
$$
$$
|||\kappa(\overline{\HH})|||^+ = \int_{\R \times S^1} \dot{\eta}(s) \Big(\max_{p \in M} H(t,p)\Big) \, \,ds\, dt = \|H\|^+ ,
$$
and
$$
|||\kappa(\overline{\HH})|||^- = -\int_{\R \times S^1} \dot{\eta}(s) \Big(\min_{p \in M} H(t,p)\Big) \,\,ds\, dt =\|H\|^-.
$$
\end{Example} 

For Hamiltonian data $(H,J)$, we will choose  a pair of homotopy triples
$$
\Hh=(\HH_L,\HH_R).
$$
This will be referred to as our {\bf cap data}.
The norm of the curvature of $\Hh$ is defined to be
$$
|||\kappa(\Hh)||| = |||\kappa(\HH_R)|||^- + |||\kappa(\HH_L)|||^+.
$$

\subsection{Floer caps}
Let $H$ be a Floer Hamiltonian.
Given a homotopy triple $\HH=(H_s,K_s,J_s)$ for $(H,J)$, we consider smooth maps 
$u \colon \R \times S^1  \to M$, which satisfy the equation
\begin{equation}\label{left-section}
    \partial_s u- X_{K_s}(u)+ J_s(u)(\partial_tu - X_{H_s}(u))=0.
\end{equation}
When $s \gg 0$ this is Floer's equation
\begin{equation*}\label{}
    \partial_s u + J(u)(\partial_tu - X_{H}(u))=0,
\end{equation*}
and for $s \ll 0$ it becomes the nonlinear Cauchy-Riemann equation
\begin{equation*}\label{}
    \partial_s u+ J^-(u)(\partial_tu)=0.
\end{equation*}
The maps satisfying \eqref{left-section} are pseudo-holomorphic sections of the bundle ${\R \times S^1} \times M \to \R \times S^1$
with respect to the almost complex structure on the total space given by
$$
 \left(
\begin{array}{ccc}
  0 & -1 & 0 \\
  1 & 0 & 0 \\
  X_{H_s}- J_s (X_{K_s}) & X_{K_s}- J_s (X_{H_s}) & J_s
\end{array}
\right).
$$

The energy of a solution $u$ of \eqref{left-section} is defined as
$$
  E(u) = \int_{\R \times S^1} \om(u) \Big( \p_su -X_{K_s}(u), J_s (\p_su -X_{K_s}(u))\Big) \,ds \, dt
$$
If the energy of $u$ is finite, then it follows from the usual arguments that
$$ u(+\infty) \eqdef \lim_{s \to \infty} u(s,t) = x(t) \in
\PP(H)$$
and
$$
u(-\infty) \eqdef \lim_{s \to -\infty} u(s,t) = p \in M.
$$
Here, convergence is in $C^{\infty}(S^1,M)$ and the point $p$ in $M$
is identified with the constant loop $t \mapsto p$. 

This asymptotic behavior implies that 
a solution $u$ of \eqref{left-section} with finite energy
determines an asymptotic spanning disk for the $1$-periodic orbit $u(+\infty)=x$ (see Figure 1). 
More precisely, for sufficiently large $s>0$, one can complete and reparameterize $u|_{[-s,s]}$ to be  
a spanning disc for $x$ in a homotopy class  which is independent of $s$. Since they play the same role for us, 
we will not distinguish between spanning discs and asymptotic spanning discs.
\begin{figure}[h]
\label{fig:leftcap} 
\caption{A left Floer cap $u$ asymptotic to $x$.}
\vspace{.5cm}
\begin{center}
\psfrag{$p$}[][][0.8]{$p$}
\psfrag{$x$}[][][0.8]{$x$}
\psfrag{$u$}[][][0.8]{$u$}
\includegraphics{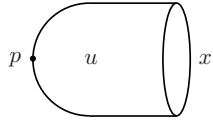}
\end{center}
\end{figure}

The set of {\bf left Floer caps} of $x \in \PP(H)$ with respect to $\HH$ is
$$\LL(x;\HH)= \left\{\ u \in C^{\infty}({\R \times S^1},M) \mid u
\text{ satisfies \eqref{left-section} },\,E(u)< \infty,\, u(+\infty)=x \right\}.
$$
For each $u \in \LL(x;\HH)$
we define the action of $x$ with
respect to $u$ by
$$
\AC_H(x,u) = \int_0^1 H(t,x(t))\,dt - \int_{\R \times S^1} u^*\om.
$$
Each left Floer cap $u \in \LL(x;\HH)$ also determines a unique homotopy class of  trivializations of $T^*M|_x$ and
hence a Conley-Zehnder index $\CZ(x,u)$.

For a generic homotopy triple $\HH$, i.e.,  a generic choice of $J_s$, each $\LL(x;\HH)$ 
 is a smooth manifold. The dimension of the component containing $u$ is $n- \CZ(x,u)$, \cite{pss}.
Hence, for every orbit $x \in \PP(H)$ with a (regular) left Floer cap $u$ we have 
$
\CZ(x,u) \leq n.
$

For any function of the form $F(s,\cdot)$, we will use the notation  $$\overleftarrow{F}(s, \cdot)= F(-s, \cdot).$$
Given a homotopy triple $\HH=(H_s,K_s,J_s)$, we will also consider maps $v \colon
{\R \times S^1} \to M$ which satisfy the equation
\begin{equation}\label{right-section}
    \partial_sv + X_{\overleftarrow{K_s}}(v)+ \overleftarrow{J_s}(v)(\partial_tv - X_{\overleftarrow{H_s}}(v))=0.
\end{equation}
In this way, we obtain for each $x \in \PP(H)$ the space of {\bf right
Floer caps},
$$\RR(x;\HH)=\left\{\ v \in C^{\infty}({\R \times S^1},M) \mid v \text{ satisfies \eqref{right-section}} ,\,E(v)< \infty,\,
v(-\infty)=x \right\}.$$ 
Every right Floer cap $v \in \RR(x;\HH)$ also determines an asymptotic spanning disc for $x$,  
$$\overleftarrow{v}(s,t) =v(-s,t).$$  
For a generic homotopy triple the space of right Floer caps $\RR(x;\HH)$,
is also a smooth manifold. The  dimension of the component containing $v$ is
 $\CZ(x,\overleftarrow{v}) -n$, and so for every orbit $x \in \PP(H)$ with a (regular) right Floer cap $v$ we 
 have 
$
\CZ(x,\overleftarrow{v}) \geq -n.
$


\subsection{Standing assumptions and their implications}

We now specify three conditions on the Hamiltonian data $(H,J)$ 
and the cap data $\Hh$ which will be assumed to hold throughout 
Section 2. We also describe some 
implications of these assumptions which will be used repeatedly.

Our first assumption, which has already been stated, is
\begin{description}
  \item[A1] $H$ is a Floer Hamiltonian.
  \end{description}
This implies that  the elements of $\PP(H)$ are isolated and hence finite
in number.  It also ensures that the limit as $s \to +\infty$ (resp. $-\infty$) 
of every left (resp. right) Floer cap is a unique element of $\PP(H)$.

Secondly, we assume 
\begin{description}
  \item[A2] The cap data $\Hh$ satisfies $|||\kappa(\Hh)|||<\hbar$.
 \end{description}
Cap data with curvature satisfying this estimate exists if $\|H\|< \hbar$ (see Example \ref{ex:linear}). More generally, 
one can find such cap data if the path generated by $H$ satisfies 
$\bar{\rho}([\phi^t_H]) < \hbar$.\footnote{This is established in the proof of Theorem \ref{technical}.} 

As described in the next subsection, the energy of the Floer caps we will consider is bounded above by $|||\kappa(\Hh)|||$
 (see inequality \eqref{uniform-energy}). 
Together with condition 
{\bf A2}, this fact will 
allow us to avoid bubbling in our compactness statements. To achieve 
this, we must first restrict our choice of the families of almost complex 
structures, $J_s$, which appear in our cap data.
We begin with the following simple observation. 

\begin{Lemma}
\label{open}
For every $\delta>0$ there is a nonempty  open subset $\JJ^{\delta} \subset \JJ(M,\om)$ such that 
for every $J \in \JJ^{\delta}$ we have $\hbar(J) \geq \hbar - \delta$.
\end{Lemma}

\begin{proof}
Assume that no such subset exists. Recall that $\JJ(M,\om)$ with its usual $C^{\infty}$-topology is a complete metric space with a metric 
defined in terms of a fixed background metric on $M$.
Choose an almost complex structure $J_0$ in $\JJ(M,\om)$ with $\hbar(J_0)\geq \hbar - \delta/2$.
Let $B(J_0,n)$ be the ball in $\JJ(M,\om)$ with radius $1/n$ and center $J_0$. By our assumption, in each open ball $B(J_0,n)$
there is a almost complex structure $J_n \in \JJ(M, \om)$ and a nonconstant $J_n$-holomorphic sphere $u_n \colon S^2 \to M$ with 
$\om(u_n) < \hbar-\delta$. By the compactness theorem for holomorphic curves (see, for example \cite{McDSa04} Theorem $4.6.1$), 
there is a subsequence of the $u_n$ which, after reparameterization, converges modulo bubbling to 
a nonconstant $J_0$-holomorphic sphere with energy no greater than $\hbar-\delta$. This contradicts the condition
$\hbar(J_0)\geq \hbar - \delta/2$.
\end{proof}

Let $\JJ^{\delta}$ be a set as in Lemma \ref{open}. Since $\JJ^{\delta}$ is open in 
the $C^{\infty}$-topology, it follows from 
\cite{fhs} that one can achieve transversality for the moduli spaces of Floer caps 
using only families  $J_s$ which take values in $\JJ^{\delta}$. In other words, these moduli spaces
can be assumed to be smooth manifolds such that the dimensions of their components agree
with their virtual dimensions.

Our final assumption is
\begin{description}
\item[A3] The families $J_s$ which appear in the cap data $\Hh$ only take values in a fixed set 
$\JJ^{\delta_{\Hh}}$, as described in Lemma \ref{open}, for
\begin{equation*}
\label{ }
\delta_{\Hh} = \frac{\hbar -|||\kappa(\Hh)|||}{2}. 
\end{equation*} 
Moreover, transversality holds for every space of left and right Floer caps for $\Hh$.
\end{description}

By condition {\bf A3}, any bubble which forms from 
a sequence of  Floer caps is a nonconstant $J$-holomorphic sphere
for some $J$ in $\JJ^{\delta_{\Hh}}$. By the definition of $\JJ^{\delta_{\Hh}}$ and 
condition {\bf A2} we have $$\hbar(J) \geq \12 (\hbar +|||\kappa(\Hh)||)> |||\kappa(\Hh)|||.$$ 
Hence, the energy of this bubble must be greater than $|||\kappa(\Hh)|||$. 
However, as mentioned above, the energy of the Floer caps 
we will consider is less than $|||\kappa(\Hh)|||$, so no such bubbling can occur. 

Since we are avoiding all holomorphic spheres in our compactifications, 
we can also fix coherent orientations on our moduli 
spaces, as in \cite{fh}. In particular, we can count zero dimensional moduli spaces with 
signs, and use the fact that the signed count of boundary components of a compact one dimensional 
moduli space is zero. 

\begin{Remark}\label{assumption order}
The first two conditions do not depend on the choice of the families of almost complex structures 
which appear in $\Hh$.  Hence, if  {\bf A1}  and {\bf A2} hold, one can always 
choose these families of almost complex structures so that {\bf A3} is satisfied. 
\end{Remark}

\subsection{Central orbits}

We say that an orbit $x \in \PP(H)$ is {\bf{central}} with respect to $\Hh$, if there is a pair 
of Floer caps $(u,v)$ in $\LL(x;\HH_L) \times
 \RR(x;\HH_R)$ such that the class  $[u \#v] \in \pi_2(M)$ is trivial.  Here, $u \# v$ denotes
the obvious concatenation of maps.
In this case, the pair $(u,v)$ will be referred to as a {\bf central pair} of Floer caps for $x$.
The point of these definitions is that there are useful bounds on the indices and actions of central 
orbits, as well as the energy of Floer caps which appear in central pairs. 

\begin{figure}[h]
\label{fig:central} \caption{A central orbit $x$ with central pair $(u,v)$.}
\vspace{.5cm}
\begin{center}
\psfrag{$v$}[][][0.8]{$v$}
\psfrag{$x$}[][][0.8]{$x$}
\psfrag{$u$}[][][0.8]{$u$}
\includegraphics{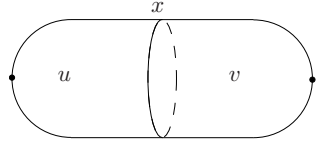}
\end{center}
\end{figure}

The dimension formulas for left and right Floer caps  imply that for
any $u \in \LL(x;\HH_L)$ and  $v \in \RR(x;\HH_R)$ we have 
$ \CZ(x,u) \leq n$ and  $\CZ(x,\overleftarrow{v}) \geq -n .$
For a central pair $(u,v)$ for $x$ we then get 
\begin{equation}
\label{index-bounds}
-n \leq \CZ(x,\overleftarrow{v}) = \CZ(x,u) \leq n.
\end{equation}

A straight forward computation  also shows that  
\begin{equation}\label{energy-left}
    0 \leq E(u) =   - \AC_H(x,u)+
\int_{\R \times S^1} \kappa(\HH_L)(s,t,u)\,ds \, dt,
\end{equation}
and
\begin{equation}\label{energy-right}
    0 \leq E(v) = \AC_H(x,\overleftarrow{v})  - \int
_{\R \times S^1} \kappa(\HH_R)(s,t,\overleftarrow{v})\,ds \,dt.
\end{equation}

If $(u,v)$ is a central pair for $x$,  then  $\AC_H(x,u) = \AC_H(x,\overleftarrow{v})$ and equations
 \eqref{energy-left} and \eqref{energy-right} imply that  
\begin{equation}\label{action-bounds}
    -|||\kappa(\HH_R)|||^- \leq \AC_H(x,\overleftarrow{v})= \AC_H(x,u) \leq  |||\kappa(\HH_L)|||^+ .
\end{equation}

For any central pair $(u,v)$ one obtains from \eqref{energy-left}, \eqref{energy-right} and \eqref{action-bounds} 
the uniform energy bounds
\begin{equation}\label{uniform-energy}
E(u), E(v) \leq |||\kappa(\Hh)|||.
\end{equation}

\subsection{A lower bound for the number of central periodic orbits}

We now prove that there are at least $\rk(\H(M;\Z))$ central periodic 
orbits of $H$ with respect to $\Hh$.  The argument relies heavily 
on the assumption,  {\bf A2},  that $|||\kappa(\Hh)|||< \hbar$.

Fix a Morse function $f$ on $M$ and a metric $g$ so that the pair $(f,g)$
is Morse-Smale. In other words, the Morse complex, 
$(\CM(f), \p_g)$, is well-defined. The chain group, $\CM(f)$,
 is the $\Z$-module generated by the
critical points of $f$, and is graded by the Morse index. The boundary map $\p_g$ counts solutions $\sigma \colon \R \to M$ of the 
ordinary differential equation
\begin{equation}\label{gradient}
\dot{\sigma}(s) = - \nabla_g f(\sigma (s)).
\end{equation}
To be more precise, for critical points $p$ and $q$, let
$$
m(p,q) =\left\{ \sigma \colon \R \to M \mid \sigma \text{ satisfies \eqref{gradient}},\, \sigma(-\infty)=p,\,
\sigma(+\infty)=q \right\}.
$$
The Morse boundary map counts  the elements  
of $m(p,q)/\R$ for $\ind(p)=\ind(q)+1$, where $\R$ acts (freely) by translation 
on the domains of the maps in $m(p,q)$.

Let $\CF(H)$ be the $\Z$-module generated by the 
elements of $\PP(H)$. To detect central orbits, we will construct  two 
$\Z$-module homomorphisms, $$\sl \colon \CM(f) \to \CF(H)$$ and  $$\sr \colon \CF(H) \to \CM(f),$$ whose
 composition $\Phi_{\Hh}= \Phi_R \circ \Phi_L$ is chain homotopic to the identity. These maps are  
similar to those constructed in Lagrangian Floer theory by Chekanov in \cite{ch}. 
In place of Floer continuation trajectories, we use the hybrid moduli spaces of \cite{pss}.

We begin by defining $\sl$. A left or right Floer cap is called {\bf short} if its energy is less than $\hbar$.
The subset of short elements in $\LL(x;\HH_L)$ is denoted by $\LL'(x;\HH_L)$.
Consider the space of left-half gradient trajectories;
$$
\ell(p) =\left\{ \alpha \colon (-\infty,0] \to M \mid \dot{\alpha} = - \nabla_g f(\alpha),\, \alpha(-\infty)=p \right\}.
$$
For a critical point $p$ of $f$ and an orbit $x$ in $\PP(H)$, set 
$$
\LL(p,x;f,\HH_L) =\{ (\alpha,u) \in \ell(p) \times \LL'(x;\HH_L) \mid \alpha(0)=u(-\infty)\}.
$$ 
For generic data, $\LL(p,x;f,\HH_L)$ is a smooth manifold and the local dimension of the 
component containing $(\alpha,u)$ is $\ind(p)-n - \CZ(x,u)$, \cite{pss}. 
The  homomorphism
$
\sl \colon \CM(f) \to \CF(H)
$ 
is now defined on each critical point $p$ of $f$ by
$$
\sl(p) = \sum_{x \in \PP(H)}  \# \LL_0(p, x ;f,\HH_L) x,
$$
where, $\# \LL_0(p, x ;f,\HH_L)$ is the number of zero-dimensional components in 
$\LL(p, x ;f,\HH_L)$ counted with signs determined by a fixed coherent orientation. 
The shortness assumption implies that $\LL_0(p, x ;f,\HH_L)$  is compact. Hence, 
$\# \LL_0(p, x ;f,\HH_L)$ is finite and $\sl$ is well-defined.

It is now easiest to define the composite homomorphism $\Phi_{\Hh}$. Consider the 
space of right-half gradient trajectories
$$
r(q) =\left\{ \beta \colon [0, +\infty) \to M \mid \dot{\beta} = - \nabla_g f(\beta),\, \beta(+\infty)=q \right\},
$$
and let $\RR'(x;\HH_R)$ be the collection of short right Floer caps of $x$.
The spaces
$$
\RR(x,q;\HH_R,f) =\{ (v, \beta) \in  \RR'(x;\HH_R)\times r(q) \mid v(+\infty) = \beta(0)\},
$$
are also smooth manifolds for generic data, and the dimension of the 
component containing $(v, \beta)$ is $\CZ(x,\overleftarrow{v}) - \ind(q)+n$. 
Let $\RR_0(x,q ; \HH_L, f)$ be the set of zero-dimensional components in $\RR(x,q ; \HH_L, f)$.
The map $\Phi_{\Hh} \colon \CM(f) \to \CM(f) $  is then defined by setting the coefficient 
of $q$ in $\Phi_{\Hh}(p)$ to be the integer
$$
 \sum_{x \in \PP(H)} \# \Big\{  ((\alpha, u), (v,\beta)) \in \LL_0(p,x;\HH_R,f)) \times \RR_0(x,q;\HH_R,f)) \mid [u\#v]=0 \Big\}.
$$ 
The map $\Phi_{\Hh}$ preserves the grading by the Morse index.
\begin{figure}[h]
\label{fig:rigid} \caption{Rigid configurations counted by $\Phi_{\Hh}$.}
\vspace{.5cm}
\begin{center}
\psfrag{alpha}[][][0.8]{$\alpha$}
\psfrag{$v$}[][][0.8]{$v$}
\psfrag{$x$}[][][0.8]{$x$}
\psfrag{$u$}[][][0.8]{$u$}
\psfrag{$p$}[][][0.8]{$p$}
\psfrag{$q$}[][][0.8]{$q$}
\psfrag{beta}[][][0.8]{$\beta$}
\includegraphics{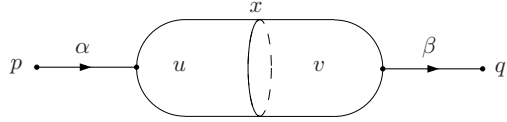}
\end{center}
\end{figure}
 
Finally, the map $\Phi_R \colon \CF(H) \to \CM(f)$ is determined by $\Phi_L$ and $\Phi_{\Hh}$
as follows.
Let $V_L$ be the submodule of $\CF(H)$ generated by the orbits $x$ in $\PP(H)$ 
which appear in an element in the image of $\sl$ with a nonzero 
integer coefficient. The maps $\sl$ and $\Phi_{\Hh}$, together with the condition
$\Phi_{\Hh} = \sr \circ \sl$, uniquely determine 
the restriction of  $\sr$ to $V_L$. 
Setting $\sr=0$ on the complement of $V_L$ we obtain the 
full map  $\sr$. In particular, the coefficient of $q$
in $\sr(x)$ is the (signed) count of elements $(v,\beta) \in   \RR_0(x,q;\HH_R,f))$ for which there is an element 
$(u, \alpha)$ in some $\LL_0(p,x;\HH_R,f))$ such that $[u\#v]=0$. 

\begin{Proposition}\label{prop:ident}
The map $\Phi_{\Hh} \colon \CM(f) \to \CM(f)$ is chain homotopic to the identity.
\end{Proposition}

\begin{proof}

We first use the cap data, $\Hh=(\HH_L,\HH_R)$, to construct a \emph{ homotopy of homotopy triples} $\{\HH^{\lambda}\} = \{(H_s^{\lambda}, K_s^{\lambda}, J_s^{\lambda})\}$ for $\lambda \in [0,+\infty)$. For convenience, we introduce the notation $\HH = \HH(s)=(H_s,K_s,J_s)$ to emphasize the $s$-dependence of $\HH$. 
Set
 $$
\HH^{1}(s) =  \begin{cases}
   \HH_{L}(s+c_1)& \text{ when $ s \leq 0$}, \\
   \HH_{R}(s-c_1)& \text{ when $ s \geq 0$},
   \end{cases}
$$
for a constant $c_1>0$ which is large enough to ensure that the domains on which $\HH_L(s+c_1)$ and $\HH_R(s-c_1)$
depend on $s$, do not intersect. For $\lambda \in [1,+\infty)$, we then define 
$$
\HH^{\lambda}(s) =  \begin{cases}
   \HH_{L}(s+c(\lambda))& \text{ when $ s \leq 0$}, \\
   \HH_{R}(s-c(\lambda))& \text{ when $ s \geq 0$},
   \end{cases}
$$
where $c(\lambda)$ is a smooth nondecreasing function which equals $\lambda$ for $\lambda \gg 1$ and is 
equal to $c_1$ for $\lambda$ near $1$.  Finally for $\lambda \in [0,1]$, we set 
$$\HH^{\lambda}(s)= (\zeta(\lambda)H^1_s, \zeta(\lambda)K^1_s, J^{\lambda}_s)$$ 
for a smooth nondecreasing function $\zeta \colon [0,1] \to [0,1]$ which equals zero near 
$\lambda = 0$ and equals one near  $\lambda=1$. As well, the compact homotopies $J^{\lambda}_s$ are chosen so that they equal $J^0_s$  for $\lambda$ near zero
and equal $J^1_s$ for $\lambda$ near one. Furthermore, the 
$J^{\lambda}_s$ are chosen so that each almost complex structure appearing in these families 
lies in $\JJ^{\delta_{\Hh}}$.

The following properties of $\HH^{\lambda}_s$ are easily verified:  
\begin{itemize}
  \item For each $\lambda\in [0,+\infty)$, $H_s^{\lambda}$ and $K_s^{\lambda}$ are compact homotopies from the zero function to itself.
  \item $\HH^0= (0, 0 , J_s^0)$.
  \item When $\lambda$ is sufficiently large $$\HH^{\lambda}(s) = \begin{cases}
   \HH_{L}(s+\lambda)& \text{for $ s \leq 0$}, \\
    \HH_{R}(s-\lambda)& \text{for $ s \geq 0$}.
   \end{cases}
$$
\item $|||\kappa(\HH^{\lambda})|||^+  = \zeta(\lambda)|||\kappa(\Hh)|||$ for all $\lambda \in [0,+\infty)$.
\end{itemize}

To construct the desired chain homotopy, we consider maps
$w$ in $\C^{\infty}({\R \times S^1},M)$ which satisfy the equation
\begin{equation}
\label{middle}
\p_s w- X_{K^{\lambda}_s}(w)+ J^{\lambda}_s(w)(\p_tw - X_{H^{\lambda}_s}(w))=0.
\end{equation}
These are perturbed holomorphic cylinders which are asymptotic, at both ends, to points in $M$.
In particular, since the perturbations are compact, each $w$ can be uniquely completed to a perturbed holomorphic sphere.
We consider these maps for all values of $\lambda \in[0,+\infty)$ and restrict our attention to those maps which represent the trivial class in $\pi_2(M)$. Let
$$\MM(\HH^{\lambda})=\left\{ (\lambda,w)  \in [0,+\infty) \times \C^{\infty}({\R \times S^1},M) \mid w \text{  satisfies  } \eqref{middle},\, [w]=0 \right\}.$$
For every $(\lambda,w) \in \MM(\HH^{\lambda})$, we have 
\begin{eqnarray*}
E(\lambda,w) & = &  \int_{\R \times S^1} \om \left( \p_s w -X_{K^{\lambda}_s}(w), J^{\lambda}_s (\p_sw -X_{K^{\lambda}_s}(w))\right) \,ds \, dt\\
& = &  \int_{\R \times S^1} \Big( \om(X_{H^{\lambda}_s}(w), \p_sw)-\om(X_{K^{\lambda}_s}(w), \p_tw)+\om(X_{K^{\lambda}_s}(w),X_{H^{\lambda}_s}(w)) \Big) \,ds \, dt\\
 & =&  \int_{\R \times S^1} \kappa(\HH^{\lambda})(s,t,w) ds \, dt \\
 &\leq& |||\kappa(\HH^{\lambda})|||^+.
 \end{eqnarray*}
The last property of $\HH^{\lambda}$ listed above, $|||\kappa(\HH^{\lambda})|||^+  = \zeta(\lambda)|||\kappa(\Hh)|||$, implies that 
for every $(\lambda,w) \in \MM(\HH^{\lambda})$ we have the uniform energy bound
\begin{equation}
\label{energy bound}
E(\lambda,w) \leq \hbar.
\end{equation}

For a pair of critical points $p$ and $q$ of $f$, consider the set $\MM(p,q;\{\HH^{\lambda}\})$ 
defined by
$$
\Big\{ \big(\alpha,(\lambda,w ), \beta \big) \in \ell(p)  \times \MM( \HH^{\lambda}) \times r(q)\mid
\alpha(0)= w(-\infty),\, w(+\infty)=\beta(0)\Big\}
$$ 
For generic data, each $\MM(p,q;\{\HH^{\lambda}\})$
is a manifold of dimension $\ind(p) - \ind(q)+1$.
We define the homomorphism $h \colon \CM(f) \to \CM(f)$ by 
$$
h(p) =\sum_{\CZ(q) = \CZ(p)+1} \# \MM(p,q;\{\HH^{\lambda}\}) q.
$$
The uniform bound \eqref{energy bound} precludes 
bubbling and implies that the zero dimensional spaces 
counted  by $h$ are compact. 

To prove that $h$ is the desired chain homotopy, it suffices 
to show that for every critical point $p$ of $f$ we have 
\begin{equation*}
\label{chain homotopy}
\big(\id - \Phi_{\Hh} + h \circ \p_g + \p_g \circ h\big)(p)=0.
\end{equation*}
Let $r$ be a critical point of $f$
with $\ind(r)=\ind(p)$. We prove that the coefficients 
of $r$ in $\big(\id - \Phi_{\Hh} + h \circ \p_g + \p_g \circ h\big)(p)$  is zero.

Consider the compactification $\overline{\MM}(p,r;\{\HH^{\lambda}\})$ 
of the one dimensional moduli space $\MM(p,r;\{\HH^{\lambda}\})$. 
The energy bound \eqref{energy bound} again prevents bubbling
for sequences in  $\MM(p,r;\{\HH^{\lambda}\})$. The
gluing and compactness theorems of Floer theory  then imply that the boundary
of $\overline{\MM}(p,r;\{\HH^{\lambda}\}) $ can be identified with the union of the  
following four compact, zero-dimensional manifolds:
\begin{enumerate}
\renewcommand{\theenumi}{\roman{enumi}}
\renewcommand{\labelenumi}{(\theenumi)}
  \item $m(p,r)$,\\ 
  \item $\bigcup_{x \in \PP(H)} \left\{ ((\alpha,u),(v,\beta)) \in \LL_0(p,x;f,\HH_L) \times \RR_0(x,r;\HH_R,f) \mid [u\#v]=0 \right\},$\\ 
  \item $\bigcup_{\CZ(q) = \CZ(p)-1} m(p,q)/\R \times \MM(q,r;\{\HH^{\lambda}\}),$\\
  \item $\bigcup_{\CZ(q) = \CZ(p)+1} \MM(p,q;\{\HH^{\lambda}\}) \times m(q,r)/\R.$
\end{enumerate}

The signed count of these boundary terms is equal to zero. 
It is clear from the definitions,
that the signed counts of the elements in sets (iii) and (iv) correspond to the 
coefficient of $r$ in $h \circ \p_g(p)$ and $\p_g \circ h(p)$, respectively.
As well, the set (i)
is empty if $p \neq r$, and  consists only of the 
constant map when $p=r$. Hence, the signed count of elements in $m(p,r)$ is equal to the 
coefficient of $r$ in $\id(p)$.

It remains to show that the signed count of elements in set (ii) is the 
coefficient of $r$ in $\Phi_{\Hh}(p)$. From the definition of $\Phi_{\Hh}$,
we only  need to show that for each tuple $((\alpha,u)(v,\beta))$ in set (ii), both $u$ and $v$
are short. Since $[u\#v]=0$, the pair $(u,v)$ is central and inequality  \eqref{uniform-energy}
yields the shortness of $u$ and $v$.
\end{proof}

\begin{Proposition}\label{prop:inner}
There are at least $\rk(\H(M;\Z))$ central periodic orbits of $H$.  In particular, there are at least $\rk(\H(M;\Z))$ contractible $1$-periodic orbits $x_j$ of $H$ which admit
spanning disks $u_j$ such that 
$$ -n\leq\CZ(x_j,u_j) \leq n$$
and
$$ -|||\kappa(\HH_R)|||^-  \leq \AC_H(x_j,u_j) \leq |||\kappa(\HH_L)|||^+ .$$

\end{Proposition}

\begin{proof}
Recall that $V_L$ is the submodule of $\CF(H)$ generated by periodic orbits $x$ which appear 
in an element in the image of $\sl$ with a nonzero coefficient. Let $K_R$ be the submodule of $\CF(H)$ generated by periodic orbits which lie in the kernel 
of $\sr$ and let $p \colon V_L \to V_L / K_r$ be the projection map. We then have
\begin{equation}
\label{break}
\Phi_{\Hh} = \sr \circ p \circ \sl.
\end{equation}
Note that any periodic orbit which appears in the image of $p \circ \sl$ is central
with respect to $\Hh$.

The homology of the Morse complex, $\H(\CM(f),\p_g)$, is isomorphic to $\H(M;\Z)$.  
Fix a basis $\{a_i\}$ for $\H(\CM(f),\p_g)$ and a set of representatives $\{A_i\}$
such that $[A_i] =a_i$.  
By Proposition \ref{prop:ident}, each
$A'_i = \Phi_{\Hh}(A_i)$ represents the same class
as $A_i$. Let $W$ (resp. $W'$) be the free submodule of rank $\rk(\H(M;\Z))$ generated
by the $A_i$ (resp. $A'_i$). 
The restriction of $\Phi_{\Hh}$ to $W$ is then an isomorphism from $W$ to $W'$.
Together with \eqref{break}, this implies that 
$$
\rk( p \circ \sl(W)) = \rk (W) = \rk(\H(M;\Z)).
$$
Hence, there are at least $\rk(\H(M;\Z))$ central periodic orbits of $H$. The bounds on the Conley-Zehnder indices
and actions follow immediately from \eqref{index-bounds} and \eqref{action-bounds} . 
\end{proof}



\subsection{Identifying Floer caps for homotopic Hamiltonian paths}

We now describe some results which allow one to choose useful capping 
data.  
 
\begin{Proposition}\label{prop:identify}
Let $\GG$ be a homotopy triple  for $(G,J_G)$. For any  $H$ in $C^{\infty}_0([\phi^t_G])$ 
there is a $J_H \in \JJ_{S^1}(M,\om)$, a homotopy triple $\widetilde{\HH}$ for $(H,J_H)$, 
and bijections
$$\PHG^{\LL} \colon \LL(x; \GG) \to \LL(\phi^t_{H} \circ (\phi^t_G)^{-1}(x);\widetilde{\HH} )$$ 
and 
$$\PHG^{\RR} \colon \RR(x; \GG) \to \RR(\phi^t_{H} \circ (\phi^t_G)^{-1}(x);\widetilde{\HH} ).$$
Moreover, these bijections preserve actions and Conley-Zehnder indices.
\end{Proposition}

\begin{proof}
Let $\varrho_t = \phi^t_{H} \circ (\phi^t_G)^{-1}$.
The map 
$\PHG \colon C^{\infty}(S^1,M) \to C^{\infty}(S^1,M)$ defined by
$$
\PHG (x)(t)= \varrho_t ( x(t)).
$$
takes contractible loops to contractible loops and hence $\PHG (\PP(G))= \PP(H).$

Set $$J_H= d \varrho_t \circ J_G \circ d (\varrho_t^{-1}).$$ 
Given  the homotopy triple 
$\GG= (G_s,K_s,J_s)$ for $(G,J_G)$, we now define 
the desired  homotopy triple $\widetilde{\HH}$ for $(H,J_H)$ and the bijection 
$$\PHG^{\LL} \colon \LL(x; \GG) \to \LL(\PHG(x);\widetilde{\HH} ).$$ 
The construction of the  bijection $\PHG^{\RR}$
is entirely similar.  
   
Let $F_s$ be a compact homotopy from $G$ to $H$ such that $F_s$ belongs to $C^{\infty}_0([\phi^t_G])$
for each $s \in \R$.
For a sufficiently large $\tau>0$, we may assume that the compact homotopies $G_s$ and $F_s$
are both independent of $s$ for $|s|>\tau$.  Thus,
$$
\widetilde{H}_s = 
\left\{
  \begin{array}{lll}
    0 & \hbox{for $s \leq -\tau$;}\\  
    G_s & \hbox{for $-\tau \leq s\leq \tau $;} \\
    F_{s-2\tau}, & \hbox{for $\tau \leq s \geq 3\tau $;}\\
    H, & \hbox{for $s \geq 3\tau $;}
   \end{array}
\right.
$$
is a compact homotopy from the zero function to $H$.

Now consider the family of contractible Hamiltonian loops 
$$\varrho_{s,t}= \phi^t_{\widetilde{H}_s} \circ (\phi^t_{G_s})^{-1}.$$ 
For each value of $s$, $\varrho_{s,t}$ is a loop based at the identity,
and 
$$
\varrho_{s,t}= 
\left\{
  \begin{array}{ll}
    \id, & \hbox{for $s \leq \tau$;} \\
    \phi^t_{F_{s-2\tau}} \circ (\phi^t_{G})^{-1}, & \hbox{for $\tau \leq s \leq 3\tau$;} \\
    \varrho_t, & \hbox{for $s \geq 3\tau$.}
   \end{array}
\right.
$$

From $\varrho_{s,t}$ we obtain two families of normalized Hamiltonian, $A_s$ and $B_s$, defined by  
$$ \p_s (\varrho_{s,t}(p)) = X_{A_s}(\varrho_{s,t}(p))$$
and
$$ \p_t (\varrho_{s,t}(p)) = X_{B_s}(\varrho_{s,t}(p)).$$
The standard composition formula for Hamiltonian flows implies that
$${B_s}= \widetilde{H}_s -G_s \circ \varrho^{-1}_{s,t},$$
where $(G_s \circ \varrho^{-1}_{s,t}) (t,p) = G_s(t,\varrho^{-1}_{s,t}(p))$. From \cite{ba}, we also 
have the following useful relation 
\begin{equation}
\label{banyaga}
\p_sB_s -\p_tA_s + \{B_s,A_s\} =0.
\end{equation}

Define
$$\widetilde{\HH} =(\widetilde{H}_s, \widetilde{K}_s ,\widetilde{J}_s )=(\widetilde{H}_s,  A_s + K_s \circ \varrho^{-1}_{s,t}, d  \varrho_{s,t} \circ J_s \circ d (\varrho_{s,t}^{-1})).$$
It is easy to verify that $\widetilde{H}$ is a homotopy triple for $(H,J_H)$.
We claim that the map $\PHG^{\LL}$ defined on $\LL(x; \GG)$ by
$$\PHG^{\LL}(u)(s,t) = \varrho_{s,t}( u(s,t))$$
is a bijection onto $\LL(\PHG(x);\widetilde{\HH}).$ It suffices to prove that for $\tilde{u} = \PHG^{\LL}(u)$ we have 
$$
\partial_s \tilde{u}- X_{\widetilde{K}_s}(\tilde{u})+\widetilde{J}_s(\tilde{u})(\partial_t \tilde{u} - X_{\widetilde{H}_s}(\tilde{u}))=0.
$$
This follows from the simple computation
\begin{eqnarray*}
  \p_s \tilde{u} + \widetilde{J}_s(\tilde{u})\p_t \tilde{u} &=& d(\varrho_{s,t}) \p_su +
  X_{A_s}(\tilde{u})
   + \widetilde{J}_s( \tilde{u})\left( d(\varrho_{s,t}) \p_tu + X_{B_s}(\tilde{u})\right) \\
  {} &=& d(\varrho_{s,t})\left( \p_su + J_s(u)\p_tu\right) + X_{A_s}(\tilde{u})+  \widetilde{J}_s(\tilde{u})X_{B_s}(\tilde{u})\\
  {} &=& d(\varrho_{s,t})\left(X_{K_s}(u)+ J_s(u)X_{G_s}(u)\right) + X_{A_s}(\tilde{u})+  \widetilde{J}_s(\tilde{u})X_{B_s}(\tilde{u}) \\
  {} &=& X_{K_s \circ \varrho_{s,t}^{-1}}(\tilde{u}) + \widetilde{J}_s(\tilde{u})\left( d(\varrho_{s,t})X_{G_s}(u) \right) + X_{A_s}(\tilde{u})+  \widetilde{J}_s(\tilde{u})X_{B_s}(\tilde{u}) \\
  {} &=& X_{K_s \circ \varrho_{s,t}^{-1}}(\tilde{u}) + \widetilde{J}_s(\tilde{u})X_{G_s \circ \varrho_{s,t}^{-1}}(\tilde{u}) + X_{A_s}(\tilde{u})+  \widetilde{J}_s(\tilde{u})X_{B_s}(\tilde{u})\\
  {} &=& X_{{\widetilde{K}_s}}(\tilde{u}) + \widetilde{J}_s(\tilde{u})X_{\widetilde{H}_s}(\tilde{u}).
\end{eqnarray*}

The additional facts that $$\CZ(x,u) = \CZ(\PHG(x), \PHG^{\LL}(u)) \text{  and  }\AC_G(x,u) = \AC_H(\PHG(x), \PHG^{\LL}(u)),$$ are well known 
from the work of Seidel in \cite{se}. (Another proof can be found in \cite{schw}.)
\end{proof}

\begin{Proposition}
\label{prop:curvature} For $\GG$ and $\widetilde{\HH}$ as in Proposition \ref{prop:identify},
$$|||\kappa(\widetilde{\HH})|||^{\pm} = |||\kappa(\GG)|||^{\pm}.$$
\end{Proposition}

\begin{proof}
The formulas $\widetilde{H}_s= B_s + G_s \circ \varrho^{-1}_{s,t}$ and $\widetilde{K}_s = A_s + K_s \circ \varrho^{-1}_{s,t}$, together with Banyaga's   formula, \eqref{banyaga},  imply that
\begin{equation}
\label{kappa}
\kappa(\widetilde{\HH}) = \kappa(\GG)\circ \varrho_{s,t}^{-1}  +\{B_s,K_s\circ \varrho_{s,t}^{-1}\} + d G_s[ \p_s( \varrho_{s,t}^{-1})]+ \{G_s\circ \varrho_{s,t}^{-1}, A_s\},
\end{equation}
where
$$
 \kappa(\GG)\circ \varrho_{s,t}^{-1}  = \p_s G_s\circ \varrho_{s,t}^{-1} - \p_t K_s \circ \varrho_{s,t}^{-1}+ \{G_s,K_s\}\circ \varrho_{s,t}^{-1}.
$$
The second term on the right side of \eqref{kappa} vanishes since $B_s =0$ for $s \in [-\tau,\tau]$ and 
$K_s=0$ for $|s|>\tau.$ We claim that the last two terms on the right side of \eqref{kappa} cancel.
For $s \leq \tau$ and $s \geq 3\tau$, both terms vanish since $\varrho_{s,t}$ is independent of $s$
in these ranges.
For $s \in [\tau, 3\lambda]$, we have 
\begin{eqnarray*}
d G_s\big( \p_s (\varrho_{s,t}^{-1})\big)+ \{G_s\circ \varrho_{s,t}^{-1}, A_s\}& = & d G\big( \p_s (\varrho_{s,t}^{-1})\big)+ \{G \circ \varrho_{s,t}^{-1}, A_s\} \\
 & = & d G\Big( \p_s (\varrho_{s,t}^{-1})+ d(\varrho_{s,t}^{-1})(X_{A_s}) \Big) \\
 & = & d G\Big( \p_s (\varrho_{s,t}^{-1})+ d(\varrho_{s,t}^{-1})\big[\p_s( \varrho_{s,t})(\varrho_{s,t}^{-1})\big] \Big)\\
 & = & 0.
\end{eqnarray*}
The final equality follows from the identity
$$
0= \p_s(\varrho_{s,t}^{-1}\circ \varrho_{s,t}(p)) =  \p_s (\varrho_{s,t}^{-1})(\varrho_{s,t}(p)) +d(\varrho_{s,t}^{-1})\big[\p_s( \varrho_{s,t}(p))\big].
$$
Hence, we have  $\kappa(\widetilde{\HH})=\kappa(\GG) \circ \varrho_{s,t}^{-1}$
which implies Proposition \ref{prop:curvature}.
\end{proof}

\section{Proof of Theorem \ref{general}}

With the machinery of the previous section in hand, the proof of Theorem \ref{general} is almost immediate.
Assume that $H$ is a Floer Hamiltonian such that $\|H\|<\hbar$. If $\phi^t_H$ does not minimize the 
positive Hofer length in its homotopy class then we can find a Hamiltonian $G$ in $C^{\infty}_0([\phi^t_H])$
such that $\|G\|^+ < \|H\|^+$.

Fix a family $J_G$ in $\JJ_{S^1}(M, \om)$ and let 
$$
J_H= d \big(\phi^t_{H} \circ (\phi^t_G)^{-1}\big)  \circ J_G \circ d\big(\phi^t_G \circ (\phi^t_{H})^{-1}\big).
$$
For the Hamiltonian data $(H,J_H)$, we consider the cap data 
$$\Hh_{\scriptscriptstyle{G}}=(\widetilde{\HH}_{\scriptscriptstyle{G}}, \overline{\HH}).$$ Here, $\overline{\HH}$  is a linear homotopy triple 
for $(H,J_H)$ as in Example \ref{ex:linear}, and $\widetilde{\HH}_{\scriptscriptstyle{G}}$ 
is the homotopy triple obtained by applying 
Proposition \ref{prop:identify} to a linear homotopy triple 
$\overline{\GG}$ for $(G,J_G)$. 

Proposition \ref{prop:curvature} implies
\begin{eqnarray*}
|||\kappa(\Hh_{\scriptscriptstyle{G}})||| & = & |||\kappa(\widetilde{\HH}_{\scriptscriptstyle{G}})|||^+ + |||\kappa(\overline{\HH})|||^- \\
 {}& = & |||\kappa(\overline{\GG})|||^+ + |||\kappa(\overline{\HH})|||^- \\
 {}& = & \|G\|^+ +\|H\|^-\\
 {}& < & \hbar.
\end{eqnarray*}
Hence, our data satisfies conditions {\bf A1} and  {\bf A2} and can therefore be chosen to satisfy {\bf A3},
 as well.  Applying Proposition \ref{prop:inner},  we obtain at least 
$\rk(\H(M;\Z))$ elements 
$x_j$ of $\PP(H)$ which admit spanning discs
$u_j$ such that 
\begin{equation*}
\label{index}
-n \leq \CZ(x_j,u_j) \leq n
\end{equation*}
and
\begin{equation*}
\label{action}
-\|H\|^-= - |||\kappa(\overline{\HH})|||^- \leq \AC_H(x_j,u_j) \leq  |||\kappa(\widetilde{\HH}_{\scriptscriptstyle{G}})|||^+ = \|G\|^+ < \|H\|^+.
\end{equation*}

\section{Proof of Theorem \ref{degenerate}}

 Let $H$ be a possibly degenerate Hamiltonian such that 
$\|H\| <\hbar$ and $\phi^t_H$ does not minimize the 
positive Hofer length in its homotopy class. 
To detect the desired $1$-periodic orbit $y$ and spanning disc $w$ such that
$$-\|H\|^- \leq \AC_H(y,w) < \|H\|^+,$$ we need to analyze 
a sequence of Floer caps. However, to obtain both the lower and upper bounds
for the action we first need a slight generalization of the machinery from 
Section \ref{floer}. 

\subsection{Preparations}   For this discussion we assume that $H$ is a Floer Hamiltonian. 
The required generalization begins
(and essentially ends) with the type of homotopy triples, $\HH=(H_s,K_s,J_s)$, we associate to the 
Hamiltonian data $(H,J)$. We now allow $H_s$ to be a compact homotopy from
any constant function $c$ to $H$.  The notions of curvature, cap data $\Hh=(\HH_L,\HH_R)$, and 
left and right Floer caps are defined in the same way as before. 

Central orbits and central pairs are also defined as in Section 2.
Only our action and energy bounds must be adjusted.
(The bounds on the Conley-Zehnder indices are unchanged.)
Equations \eqref{energy-left}
and \eqref{energy-right}, become 
\begin{equation}\label{energy-left-new}
    0 \leq E(u) =  c_L - \AC_H(x,u)+
\int_{\R \times S^1} \kappa(\HH_L)(s,t,u)\,ds \, dt,
\end{equation}
and
\begin{equation}\label{energy-right-new}
    0 \leq E(v) = \AC_H(x,\overleftarrow{v}) - c_R - \int
_{\R \times S^1} \kappa(\HH_R)(s,t,\overleftarrow{v})\,ds \,dt,
\end{equation}
where $c_L$ and $c_R$ are the constant functions coming from $\HH_L$ and $\HH_R$, respectively.

If $x$ is central for $\Hh$ and  $(u,v)$ is a central pair for $x$, we now have action bounds  \begin{equation}\label{action-bounds-new}
    -|||\kappa(\HH_R)|||^- + c_R \leq \AC_H(x,\overleftarrow{v})= \AC_H(x,u) \leq  |||\kappa(\HH_L)|||^+ +c_L.
\end{equation}

For a central pair $(u,v)$, the previous uniform energy bound \eqref{uniform-energy} becomes
\begin{equation}\label{uniform-energy-new}
E(u), E(v) \leq |||\kappa(\Hh)||| + c_L -c_R .
\end{equation} 

If we replace assumption {\bf A2} with the assumption that
$$ |||\kappa(\Hh)||| + c_L -c_R < \hbar,$$ 
then Propositions  \ref{prop:ident} and \ref{prop:inner} hold, for the adjusted action bounds
coming from \eqref{action-bounds-new}.
In particular,  for every Morse-Smale pair $(f,g)$
one can construct a $\Z$-module homomorphism
$$
\Phi_{\Hh} \colon \CM(f) \to \CM(f)
$$
which is chain homotopic to the identity, and this implies the existence of
at least $\rk(\H(M;\Z))$ orbits in $\PP(H)$ which are central with respect to $\Hh$.

The point of this modification is that we can now use {\bf decreasing 
linear homotopy triples} of the form 
$$
\underline{\HH} = \big( -(1-\eta(s)) \|H\|^- + \eta(s)H, 0, J_s\big).
$$
Here, $\eta$ is the same function as in Example \ref{ex:linear}, and $c=-\|H\|^-$.
The curvature norms for such triples are 
\begin{equation}\label{curvature-linear}
|||\kappa(\underline{\HH})|||^+ = \|H\| 
\text{   and    }
|||\kappa(\underline{\HH})|||^- = 0.
\end{equation}

Our proof of Theorem \ref{degenerate} will rely on a special property
acquired by Floer caps defined using deceasing linear homotopy 
triples (see Lemma \ref{nonincreasing}).

\subsection{The proof}
Returning to the setting of Theorem \ref{degenerate}, 
we have $\|H\|< \hbar$ and $\|H\|^+> \rho^+([\phi^t_H])$.
Choose a  
$G$ in $C^{\infty}_0([\phi^t_H])$ such that $\|G\|^+ < \|H\|^+$.
For some $\eps>0$ we then have
$$
\|G\|^+ < \|H\|^+ -2\epsilon.
$$

Fix a family of  almost complex structures $J_H$ in $\JJ_{S^1}(M,\om)$. To apply the Floer theoretic tools
developed above, we consider a sequence of Hamiltonian data
$(H_k,J_{k})$ such that each $H_k$ is a Floer Hamiltonian and the sequence $(H_k,J_{k})$ converges to $(H,J_H)$ in the $C^{\infty}$-topology.

We may assume that $\|H_k\| < \hbar$ for all $k$. We may also assume that for each $k$
there is a Hamiltonian $G_k$ in $C^{\infty}_0[ \phi^t_{G_k}]$  such that $\|G_k\|^+ \leq \|H_k\|^+-\eps.$
To see this, consider the Hamiltonian loop $\phi^t_{H_k} \circ (\phi^t_{H})^{-1}$. It is generated 
by the Hamiltonian 
$$
F_k =  H_k - H\circ \phi^t_{H} \circ (\phi^t_{H_k})^{-1}
$$
which converges to zero in the $C^{\infty}$-topology. 
The path 
$$
\phi^t_{H_k} \circ (\phi^t_{H})^{-1} \circ \phi^t_G
$$
is homotopic to $\phi^t_{H_k}$, relative its endpoints, and is generated by the function
$$
G_k = F_k +G \circ (\phi^t_{F_k})^{-1}.
$$
The functions $G_k$ clearly converge to $G$ and hence satisfy 
\begin{equation}\label{gk hk}
|G_k\|^+ \leq \|H_k\|^+-\eps
\end{equation} 
for sufficiently large $k$, as desired.

We now proceed as in the proof of Theorem \ref{general}. 
Let $\underline{\HH_k}$ be a decreasing linear homotopy triple for $(H_k, J_k)$. Choose the families of almost complex structures $J_{k,s}$ in the $\underline{\HH_k}$ so that they converge to $J_s$ in $\JJ_{S^1}(M, \om)$. Then $\underline{\HH_k}$ converges to $$\underline{\HH}= \big( -(1-\eta(s)) \|H\|^- + \eta(s)H, 0, J_s\big)$$ in the $C^{\infty}$-topology.
By \eqref{energy-right-new} and \eqref{curvature-linear}, the action of any $x$ with respect to a right Floer cap $v \in  \RR(x;\underline{\HH_k})$
satisfies
\begin{equation}\label{new right}
\AC_{H_k}(x,\overleftarrow{v}) \geq -|||\kappa(\underline{\HH_k})|||^- -\|H_k\|^- =-\|H_k\|^-. 
\end{equation}

Let $\underline{\GG_k}$ be a decreasing 
linear homotopy triple for $G_k$. As in Proposition \ref{prop:identify}, one can construct from $\underline{\GG_k}$ a homotopy
triple $\widetilde{\HH}_{\scriptscriptstyle{{G_k}}}$ for $(H_k,J_k)$ and an action-preserving bijection between the 
corresponding spaces 
of left Floer caps. Hence, by \eqref{energy-left-new} and \eqref{curvature-linear}, the action of any $x$ with respect to any left Floer cap $u \in  \LL(x;\widetilde{\HH}_{\scriptscriptstyle{{G_k}}})$
satisfies
\begin{equation}\label{new left}
\AC_{H_k}(x,u) \leq |||\kappa(\underline{\GG_k})|||^+ -\|G_k\|^- =\|G_k\|^+. 
\end{equation}

For each $k$ we then consider the homotopy data
$$\Hh_k= (\widetilde{\HH}_{\scriptscriptstyle{{G_k}}},  \underline{\HH_k}).$$
Arguing as in Proposition \ref{prop:curvature}, one 
can show that 
$\underline{\GG_k}$ and $\widetilde{\HH}_{\scriptscriptstyle{{G_k}}}$ have the same curvature norms. Hence,
\begin{eqnarray*}
|||\kappa(\Hh_k)||| +c_{L,k} -c_{R,k} & = & ||| \kappa(\widetilde{\HH}_{\scriptscriptstyle{{G_k}}})|||^+ + |||\kappa(\underline{\HH_k})|||^- -\|G_k\|^-  + \|H_k\|^-\\
{} & = & |||\kappa(\underline{\GG_k})|||^+  -\|G_k\|^-  + \|H_k\|^-\\
{} & = & \|G_k\|^+ +\|H_k\|^-\\
{} &<& \|H_k\|\\
{} &<& \hbar.
\end{eqnarray*}
For a Morse-Smale pair $(f,g)$ 
we can then construct, for each $k$, a chain map 
$$\Phi_{\Hh_k}= \srl \circ p_k \circ \sll \colon \CM(f) \to \CM(f)$$
which is chain homotopic to the identity. This detects 
at least $rk(\H(M;\Z))$ central orbits with respect to $\Hh_k$. 
If $x$ is one of these orbits with a central pair $(u,v)$, then \eqref{gk hk},\eqref{new right} and \eqref{new left} imply that   
\begin{equation}\label{away}
-\|H_k\|^- \leq \AC_{H_k}(x,\overleftarrow{v}) = \AC_{H_k}(x,u)\leq \|H_k\|^+ -  \eps.
\end{equation}

For simplicity, let us choose the Morse function $f$  so that it
has a unique local (and hence global) maximum at $Q \in M$.
Standard arguments imply that $Q$ is a nonexact cycle of degree $2n$ in the 
Morse complex $(\CM_*(f), \p_{g})$, and that $\Phi_{\Hh_k}(Q) =Q$.
Let $$X_k = p_k \circ \sll(Q).$$  By the construction of $\Phi_{\Hh_k}$, $X_k$ is a finite sum of the form $$X_k=\sum n^j_k x^j_k$$
where the $n^j_k$ are nonzero integers and the $x^j_k$ are central 
periodic orbits of $H_k$ with respect to $\Hh_k.$ 
Since $X_k$
gets mapped to $Q$ under $\srl$, the moduli space  
$$
\RR_0(X_k,Q;\underline{\HH_k},f) = \bigcup_j\RR_0(x^j_k,Q;\underline{\HH_k},f),
$$
which determine $\srl(X_k)$, must be nonempty.  Choose a 
$(v_k, \sigma_k)$ in $\RR_0(X_k,Q;\underline{\HH_k},f)$  for each $k$. The caps $v_k$ belongs to $\RR(x^j_k;\underline{\HH_k})$  for some $x^j_k$ in $\PP(H)$ which appears in  $X_k$ with a nonzero coefficient. 
Moreover, each $v_k$ is part of a central pair for $x^j_k$, and so by \eqref{away} we have 
\begin{equation}\label{gaps}
-\|H_k\|^- \leq \AC_{H_k}(x^j_k,\overleftarrow{v}_k)\leq \|H_k\|^+-  \eps.
\end{equation}
By \eqref{uniform-energy-new} and the computation above, we also have the uniform  bound 
\begin{equation}\label{energy}
E(v_k) \leq |||\kappa(\Hh_k)||| +c_{L,k}-c_{R,k}  < \hbar.
\end{equation}

Since the  $\underline{\HH_k}$ converge to $\underline{\HH}$, 
the energy bound \eqref{energy} implies that there is a subsequence of
the $v_k$ (which we still denote by $v_k$) that converges to a solution $v \colon \R \times S^1 \to M$ of
$$
\partial_sv + \overleftarrow{J_s}(v)(\partial_tv - X_{\overleftarrow{H_s}}(v))=0,
$$
for ${\overleftarrow{H_s}} = -(1 -\eta(-s))\|H\|^- + \eta(-s)H$.
The limit $v$ also satisfies the energy bound $E(v)< \hbar$. Hence,  
$\lim_{s \to +\infty}v(s,t)$  exists and is equal to some point $P$ in $M$. 
There is also a sequence
$s_j \to -\infty$ such that $v(s_j,t)$ converges to some $y(t) \in \PP(H).$
\footnote{The negative limit of $v$ may not exist since the orbits of $H$ may be degenerate.}

For $s \in \R$, we set $$y^{\scriptscriptstyle{[s]}}(t) = v(s,t)$$ and  
$$
\overleftarrow{v}^{\scriptscriptstyle{[s]}} = \overleftarrow{v}|_{\scriptscriptstyle{(-\infty,-s]}}.
$$
Note that $\overleftarrow{v}^{\scriptscriptstyle{[s]}}$ is a spanning disc for $y^{\scriptscriptstyle{[s]}}.$
In this notation, $y^{\scriptscriptstyle{[s_j]}} \to y$ as $j \to \infty$, and for large values of $j$ the discs 
$\overleftarrow{v}^{\scriptscriptstyle{[s_j]}}$ can be easily extended to  a spanning disc $w$
for $y \in \PP(H)$ in a fixed homotopy class. It remains to show that the action of $y$ with respect to $w$ is less than 
$\|H\|^+$ and greater than or equal to $-\|H\|^-$. Since the extensions of the $\overleftarrow{v}^{\scriptscriptstyle{[s_j]}}$ can be made arbitrarily small for sufficiently large $j$, 
it suffices to prove that for all $j$ we have  
\begin{equation}\label{need}
-\|H\|^- \leq \AC_H(y^{\scriptscriptstyle{[s_j]}},\overleftarrow{v}^{\scriptscriptstyle{[s_j]}} ) \leq  \|H\|^+-  \eps.
\end{equation}

By definition,
$$y^{\scriptscriptstyle{[s_j]}}(t) = \lim_{k \to \infty} v_k(s_j,t)= \lim_{k \to \infty}y_k^{\scriptscriptstyle{[s_j]}}(t).$$
Hence, 
\begin{equation}
\label{k-limit}
\AC_H(y^{\scriptscriptstyle{[s_j]}},\overleftarrow{v}^{\scriptscriptstyle{[s_j]}} ) = \lim_{k \to \infty} 
\AC_{H_k}(y_k^{\scriptscriptstyle{[s_j]}},\overleftarrow{v}_k^{\scriptscriptstyle{[s_j]}} ).
\end{equation}

\begin{Lemma}\label{nonincreasing}
The function $s \mapsto \AC_{H_k}(y_k^{\scriptscriptstyle{[s]}},\overleftarrow{v}_k^{\scriptscriptstyle{[s]}} )$ is nonincreasing.
\end{Lemma}

\begin{proof}
As a right Floer cap for $\underline{\HH_k}$, $v_k$ is a solution of the 
equation 
$$
\partial_sv_k + \overleftarrow{J_{k,s}}(v_k)(\partial_tv_k - X_{\overleftarrow{H_{k,s}}}(v_k))=0
$$
for $$\overleftarrow{H_{k,s}} = -(1 -\eta(-s))\|H_k\|^- + \eta(-s)H_k.$$
Let $b>a$. Then
\begin{eqnarray*}
  0 &\leq& \int_0^1 \int_{a}^{b} \om(\partial_s v_k, \overleftarrow{J_{k,s}}(v_k) (\partial_s v_k)) \,ds\,dt   \\
  {} &=& \int_0^1 \int_{a}^{b} \om ( \partial_s v_k, \partial_t v_k) - \om (\partial_s v_k, X_{\overleftarrow{H_{k,s}}})\,ds\,dt  \\
  {} &=&  \AC(y_k^{\scriptscriptstyle{[b]}},\overleftarrow{v_k}^{\scriptscriptstyle{[b]}} )-\AC(y_k^{\scriptscriptstyle{[a]}},\overleftarrow{v_k}^{\scriptscriptstyle{[a]}})  - \int_0^1 \int_{a}^{b}d \overleftarrow{H_{k,s}}(\partial_s v_k) \,ds\,dt \\
  {} &=&  \AC(y_k^{\scriptscriptstyle{[b]}},\overleftarrow{v_k}^{\scriptscriptstyle{[b]}} )-\AC(y_k^{\scriptscriptstyle{[a]}},\overleftarrow{v_k}^{\scriptscriptstyle{[a]}}) - \int_0^1 \int_{a}^{b} \left( \p_s( \overleftarrow{H_{k,s}}(v_k))- \partial_s \overleftarrow{H_{k,s}}(v_k) \right) \,ds\,dt \\
  {} &=&  \AC_{H_k}(y_k^{\scriptscriptstyle{[a]}},\overleftarrow{v_k}^{\scriptscriptstyle{[a]}}) - \AC_{H_k}(y_k^{\scriptscriptstyle{[b]}},\overleftarrow{v_k}^{\scriptscriptstyle{[b]}} )   +  \int_0^1 \int_{a}^{b} \partial_s \overleftarrow{H_{k,s}}(v_k)\,ds\,dt.
 \end{eqnarray*}
Since 
$$
\partial_s \overleftarrow{H_{k,s}} = \dot{\eta}(-s)(-\|H_k\|^- -H_k) \leq 0
$$
we have 
$$
\AC_{H_k}(y_k^{\scriptscriptstyle{[a]}},\overleftarrow{v_k}^{\scriptscriptstyle{[a]}}) \geq \AC_{H_k}(y_k^{\scriptscriptstyle{[b]}},\overleftarrow{v_k}^{\scriptscriptstyle{[b]}} ).
$$
\end{proof}

Lemma \ref{nonincreasing}, together with \eqref{gaps}, implies that
\begin{eqnarray*}
  \AC_{H_k}(y_k^{\scriptscriptstyle{[s_j]}},\overleftarrow{v}_k^{\scriptscriptstyle{[s_j]}} ) &\leq& \lim_{j \to \infty} \AC_{H_k}(y_k^{\scriptscriptstyle{[s_j]}},\overleftarrow{v_k}^{\scriptscriptstyle{[s_j]}} ) \\
   {}&=& \AC_{H_k}(x^j_k,\overleftarrow{v_k}) \\
   {}&\leq& \|H_k\|^+ - \epsilon .
\end{eqnarray*}

Hence, by \eqref{k-limit}, we have 
\begin{eqnarray*}
  \AC_H(y^{\scriptscriptstyle{[s_j]}},\overleftarrow{v}^{\scriptscriptstyle{[s_j]}} ) &=&  \lim_{k \to \infty} 
\AC_{H_k}(y_k^{\scriptscriptstyle{[s_j]}},\overleftarrow{v_k}^{\scriptscriptstyle{[s_j]}} )\\
   {}&\leq & \|H\|^+-\eps.
\end{eqnarray*}

Similarly, Lemma \ref{nonincreasing} implies
\begin{eqnarray*}
  \AC_{H_k}(y_k^{\scriptscriptstyle{[s_j]}},\overleftarrow{v_k}^{\scriptscriptstyle{[s_j]}} ) &\geq& \lim_{s \to +\infty} \AC_{H_k}(y_k^{\scriptscriptstyle{[s]}},\overleftarrow{v_k}^{\scriptscriptstyle{[s]}} )\\
   {}&=&  \int_0^1 H_k(t,P) \, dt \\
   {}&\geq& -\|H_k\|^-
\end{eqnarray*}
and so
\begin{eqnarray*}
  \AC_H(y^{\scriptscriptstyle{[s_j]}},\overleftarrow{v}^{\scriptscriptstyle{[s_j]}} ) &=&  \lim_{k \to \infty} 
\AC_{H_k}(y_k^{\scriptscriptstyle{[s_j]}},\overleftarrow{v_k}^{\scriptscriptstyle{[s_j]}} )\\
   {}&\geq &- \|H\|^-.
\end{eqnarray*}
This proves  \eqref{need} and the proof of Theorem \ref{degenerate} is complete.

\section{Proof of Theorem \ref{quasi-autonomous}}

Let $H$ be a quasi-autonomous Floer Hamiltonian
on a spherically rational symplectic manifold $(M,\om)$, such that $\|H\| < r(M,\om)$. 
We also assume  that $\phi^t_H$ does not minimize $\rho_{\scriptscriptstyle{\H}}$ in its homotopy class.
Hence,  $\phi^t_H$ does not minimize at least one of the one-sided seminorms. 
We consider the case when it does not minimize $\rho^+$.
The other case can be dealt with in similar way (see Corollary \ref{flip}).

By Theorem \ref{general}, there are at least $\rk(\H(M;\Z))$ contractible $1$-periodic orbits $x_j$ of $H$
which admit spanning disks $u_j$ such that 
the action values $\AC_H(x_j,u_j) $ lie in 
the interval $[-\|H\|^-,\|H\|^+).$  Let $P$ be a fixed global maximum of $H$. It is a 
constant $1$-periodic orbit of $H$ since $dH_t(P)=0$ for all $t \in [0,1]$. We claim that $P$ is not 
one of the orbits detected by Theorem \ref{general}.

Any spanning disk $u$ for $P$ represents an element 
 $[u]$ in $\pi_2(M)$ and so 
 \begin{eqnarray*}
  \AC_H(P,u)  &=& \int_0^1H(P)\, dt - \int_{D^2} u^* \om \\
   {}&=& \|H\|^+ - \om([u]).
\end{eqnarray*}
If $ \om ([u]) \leq 0$, then have $\AC_H(P,u) \geq \|H\|^+$. Otherwise we have $\om([u]) \geq  r(M, \om ) > \|H\|^+ +  \|H\|^-$ and  $\AC_H(P,u) < -\|H\|^-$. In either case, $P$ does not admit a spanning disc with action
 in $[-\|H\|^-, \|H\|^+)$ and so does not contribute to the count of orbits detected in Theorem \ref{general}. 
Hence, there are at least $\rk(\H(M;\Z))+1 \geq \SB(M)+1$ elements of $\PP(H)$ 

The proofs of the Arnold conjecture for compact symplectic manifolds in \cite{fl1,fl2,hs,fo,lt},
imply that the number of elements in $\PP(H)$ is at least $\SB(M)$. Moreover, they imply that the 
cardinality of $\PP(H)$ is equal to $\SB(M)$ modulo $2$. In particular,  the various versions of 
Floer homology in these works are shown to be isomorphic to the Morse homology of $M$, 
with suitable coefficient rings. 
Hence, the number of generators of the Floer complex has the same parity as the number of critical points
of a Morse function. By the strong Morse inequalities this parity is equal to the parity of $\SB(M)$.
Therefore, since we have detected at least $\SB(M) +1$ orbits in $\PP(H)$, there must be at least
$\SB(M) +2$.

\section{Proof of Theorem \ref{technical}}

Suppose that $H$ is Floer Hamiltonian and that $\bar{\rho}([\phi^t_H]) < \hbar$. We may assume that $\epsilon^{\pm}>0$ satisfy $$\epsilon^+ + \epsilon^- <
\hbar -\bar{\rho}([\phi^t_H]).$$ Choose Hamiltonians 
$G$ and $F$ in $C^{\infty}_0([\phi^t_H])$ such that 
\begin{equation*}
\label{ }
\|G\|^+ < \rho^+([\phi^t_H])+\epsilon^+,
\end{equation*}
and
\begin{equation*}
\label{ }
\|F\|^- < \rho^-([\phi^t_H])+ \epsilon^-.
\end{equation*}
Fix a family $J_G$ in $\JJ_{S^1}(M, \om)$ and set 
$$
J_H= d (\phg)  \circ J_G \circ d((\phg)^{-1})
$$
and 
$$
J_F= d (\pfg)  \circ J_G \circ d((\pfg)^{-1}),
$$
for
$\phg = \phi^t_{H} \circ (\phi^t_G)^{-1}$ and $\pfg=\phi^t_{F} \circ (\phi^t_G)^{-1}$.
We then consider the cap data $\widetilde{\Hh}_{\scriptscriptstyle{GF}}=(\widetilde{\HH}_{\scriptscriptstyle{G}}, \widetilde{\HH}_{\scriptscriptstyle{F}})$ for the Hamiltonian data $(H,J_H)$,
where $\widetilde{\HH}_{\scriptscriptstyle{G}}$ 
and $\widetilde{\HH}_{\scriptscriptstyle{F}}$ are the homotopy triples obtained by applying 
Proposition \ref{prop:identify} to  linear homotopy triples 
$\overline{\GG}$ and $\overline{\FF}$ for $(G,J_G)$ and $(F,J_F)$, respectively. 
Proposition \ref{prop:curvature}, together with our choice of $G$ and $F$, implies that
\begin{eqnarray*}
|||\kappa(\widetilde{\Hh}_{\scriptscriptstyle{GF}})||| & = & |||\kappa(\widetilde{\HH}^G)|||^+ + |||\kappa(\widetilde{\HH}^F)|||^- \\
 {}& = & |||\kappa(\overline{\GG})|||^+ + |||\kappa(\overline{\FF})|||^- \\
 {}& = & \|G\|^+ +\|F\|^-\\
 {}& < & \hbar.
\end{eqnarray*}
Hence, we can again apply Proposition \ref{prop:inner} to obtain at least $\rk(\H(M;\Z))$ elements 
$x_j$ of $\PP(H)$ which admit spanning discs
$u_j$ such that 
\begin{equation*}
\label{index}
-n \leq \CZ(x_j,u_j) \leq n
\end{equation*}
and
\begin{equation*}
\label{action}
-\rho^-([\phi^t_H])-\epsilon^- \leq -\|F\|^-\leq \AC_H(x_j,u_j) \leq \|G\|^+ \leq \rho^+([\phi^t_H]) +\epsilon^+.
\end{equation*}

\begin{Corollary}[Theorem \ref{general} for $\rho^-$]\label{flip}
Suppose $H$ is a Floer Hamiltonian such that $\|H\|< \hbar$ and $\phi^t_H$ does not minimize 
the negative Hofer length in its homotopy class. There are 
at least $\rk(\H(M;\Z))$ elements 
$x_j$ of $\PP(H)$ which admit spanning discs
$u_j$ such that 
\begin{equation*}
\label{index}
-n \leq \CZ(x_j,u_j) \leq n
\end{equation*}
and
\begin{equation*}
\label{action}
- \|H\|^- < \AC_H(x_j,u_j) \leq \|H\|^+.
\end{equation*}
\end{Corollary}

\begin{proof}

Let $F$ be a Hamiltonian in $C^{\infty}_0([\phi^t_H])$
such that $\|F\|^- < \|H\|^-$. Choose $\epsilon^{\pm}=0$ and $G=H$, and repeat the proof of Theorem \ref{technical}.
\end{proof}

\end{document}